\documentclass[12pt]{article}
\RequirePackage[OT1]{fontenc}
\RequirePackage{amsthm,amsmath, amsfonts,latexsym,amssymb}
\RequirePackage[numbers]{natbib}
 \usepackage[latin1]{inputenc}

 \usepackage{amsmath,amsfonts,latexsym,amssymb,comment}
 \usepackage[symbol]{footmisc}
 \usepackage[active]{srcltx}
 \usepackage{psfrag, epsfig, graphicx}
\usepackage[frenchb]{babel}
\usepackage{lmodern}
\theoremstyle{plain}

\renewcommand{\(}{$\,}
\renewcommand{\)}{\,$}
\newcommand{\cc}[1]{\mathcal{#1}}
\newcommand{\bb}[1]{\boldsymbol{#1}}

\renewcommand{\hat}[1]{\widehat{#1}}
\renewcommand{\tilde}[1]{\widetilde{#1}}
\newtheorem{theorem}{Theorem}
\newtheorem{proposition}{Proposition}
\newtheorem{definition}{Definition}

\def\eqdef{\stackrel{\operatorname{def}}{=}}

\def\T{\top}

\def\dd{\mathrm{d}}

\def\RR{\mathbb{R}}

\def\SSd{\mathbb{S}^{d-1}}
\def\SS{\mathbb{S}^{1}}

\def\kappa{\varkappa}
\def\ta{\theta}
\def\eps{\varepsilon}

\def\EE{\mathbb{E}}

\def\TH{\operatorname{TH}}

\newcommand{\expect}[1]{\EE^{\eps}_{#1}}

\def\bob{\boldsymbol{b}}

\def\C(S){1}                                        
\def\K{\cc K}                                       

\begin{document}


\begin{center}
{\Large
	{\sc Estimation adaptative dans le mod\`{e}le \\single-index par l'approche d'oracle\footnote{Funding of the ANR-07-BLAN-0234 is acknowledged. The second author is also supported by the DFG FOR 916.}

}
}
\bigskip

Oleg Lepski$^{1}$ \& Nora Serdyukova$^{2}$
\bigskip

{\it
$^{1}$ Laboratoire d'Analyse, Topologie, Probabilit\'es UMR 7353,
Aix-Marseille Universit\'e \\39, rue F. Joliot Curie 13453 Marseille FRANCE. \\E-mail : Oleg.Lepski@cmi.univ-mrs.fr

$^{2}$ Institute for Mathematical Stochastics, Georg-August-Universit\"at G\"ottingen \\Goldschmidtstra{\ss}e 7, 37077 G\"ottingen GERMANY. \\E-mail : Nora.Serdyukova@gmail.com
}
\end{center}
\bigskip


{\bf R\'esum\'e.}
Dans le cadre de l'estimation non paramétrique
d'une fonction multidimensionnelle nous nous intéressons à
l'adaptation structurelle. Nous supposons que la fonction à estimer
possède la structure « single-index » dans laquelle ni fonction de lien ni vecteur d'indice
ne sont connus. Nous proposons une nouvelle procédure qui s'adapte
simultanément à l'indice inconnu ainsi qu'à la régularité de la
fonction de lien. Nous présentons une inégalité d'oracle « locale »
(définie par la semi-norme ponctuelle ) pour la procédure proposée,
qui est ensuite utilisée pour obtenir la borne supérieure du risque
maximal sous une hypothèse de régularité sur la fonction de lien.
D'après la borne inférieure obtenue pour le risque minimax
l'estimateur construit est un estimateur adaptatif optimal sur
l'ensemble de classes considérées. Pour la même procédure on établit
également une inégalité d'oracle « globale » (en norme $L_r$, $r< \infty $) et étudie
sa performance sur les classes de Nikol'skii. Cette étude montre que
la méthode proposée peut être appliquée à l'estimation de fonctions
ayant une régularité inhomogène.
\smallskip

\noindent {\bf Mots-cl\'es.} Estimation adaptative,  Borne inf\'{e}rieure, Vitesse minimax, In\'{e}galit\'{e} d'oracle,  Mod\`{e}le single-index, Adaptation structurelle, R\'{e}gularit\'{e} inhomog\`{e}ne. 
\bigskip

{\bf Abstract.}
In the framework of nonparametric multivariate function estimation we are interested in structural adaptation. We assume that the function to be estimated possesses the ``single-index'' structure where neither the link function nor the index vector is known. We propose a novel procedure that adapts simultaneously to the unknown index and smoothness of link function. For the proposed procedure, we present a ``local'' oracle inequality (described by the pointwise seminorm), which is then used to obtain the upper bound on the maximal risk under regularity assumption on the link function.
The lower bound on the minimax risk shows that the constructed estimator is optimally rate adaptive over the considered range of classes. For the same procedure we also establish a ``global'' oracle inequality  (under the $ L_r $ norm, $r< \infty $) and study its performance over the Nikol'skii classes. This study shows that the proposed method can be applied to estimating functions of inhomogeneous smoothness.
\smallskip




\bigskip



\par \noindent  {\bf  \emph{Model and set-up.}}
We observe  a path \( \{ Y_{\eps}(t), t \in \cc D \} \) satisfying the equation
\begin{equation}\label{WGN_model}
    Y_{\eps}(\dd t) = F(t) \dd t + \eps W(\dd t) \; , \;\; t = (t_1, \ldots, t_d) \in [-1, 1]^d,
\end{equation}
where \( W \) is a Brownian sheet and \( \eps \in (0,1) \). 
We consider $d=2$ except the second assertion of Theorem \ref{th:pointwise-adaptation} concerning a lower bound for function estimation at a
point. 
Additionally, we assume that the function \( F \) has the single-index structure, i.e. there exist an unknown link function
\( f : \RR \to \RR \) and an index vector \( \ta^* \in \SS \) such that
\begin{equation}\label{single-index}
    F(x) = f(x^{\T} \ta^*).
\end{equation}
We suppose that $f\in\mathbb{F}_M=\left\{g: \RR \to \RR\; |\; \sup_{u\in\RR}|g(u)| \le M\right\}$ for some \( M>0 \),
however its knowledge is not required for the estimation procedure.
Our aim is to estimate the entire function \( F \) on \( [-1/2, 1/2]^2 \) or its value \( F(x) \) from the observation \( \{ Y_{\eps}(t), t \in \cc D \} \) 
without any prior knowledge of the nuisance parameters \( f\) and \( \ta^{*} \).
The quality of estimation is measured by
\( \cc R_r^{(\eps)} (\hat F, F)   =    \expect{F} \| \hat F  - F \|_r \), where \( \| \cdot \|_r \) is the \( L_r \) norm on \( [-1/2, 1/2]^2 \), \( r \in [1, \infty) \), or by the ``pointwise''  risk
\( \cc R_{r,x}^{(\eps)} (\hat F, F)        =    (\expect{F} |\hat F (x) - F(x)|^r )^{1/r} \).

\medskip

\par \noindent {\bf \emph{Objectives.}}
The goal of our study is at least threefold. First, we seek an estimation procedure $\hat{F}(x), x\in[-1/2,1/2]^{2},$ for $F$ which could
be applicable to any function $F$ satisfying (\ref{single-index}).
Moreover, we want to bound the risk of this estimator uniformly
over the set $\mathbb{F}_M\times \SS $. More precisely, we establish for  $\hat{F}(x)$ the local oracle inequality:
\begin{equation}
\label{eq:local-oracle-intro}
  \cc R_{r,x}^{(\eps)} (\hat F, F)\leq C_r A^{(\eps)}_{f,\theta^*}(x),\quad\forall f\in\mathbb{F}_M,\;\;\forall\theta^*\in\SS ,\;\;\forall x\in[-1/2,1/2]^{2}.
\end{equation}
Here the quantity $A^{(\eps)}_{f,\theta^*}$ is completely determined by the function $f$, vector~$\theta^*$ and noise level $\eps$, while $C_r$ is a numerical constant independent of \( F \) and \( \eps \). Next, we apply this result to minimax adaptive estimation over the scale of $\mathbb{H}(\beta,L)$, H\"older classes of functions, see Definition \ref{def:holder-class}.
In particular, we find the minimax rate over $\mathbb{H}(\beta,L)\times \SS $ and
prove that our estimator $\hat{F}$ achieves that rate, i.e. is optimally adaptive. This result is quite surprising because, if $\theta^*$
is fixed, say $\theta^*=(1,0)^{\T}$, then it is well known that an optimally adaptive estimator does not exist, see \cite{Lep1990}.


Note also that local oracle inequality (\ref{eq:local-oracle-intro}) allows us to bound from above the ``global'' risk as well: \( \cc R_r^{(\eps)} (\hat F, F)\leq C_r \|A^{(\eps)}_{f,\theta^*}\|_r \).
The latter is a global oracle inequality.  As local oracle inequality (\ref{eq:local-oracle-intro})
is a powerful tool for deriving minimax adaptive results in pointwise estimation, so global oracle inequality
can be used for constructing adaptive estimators of the entire function \( F \). We will consider the collection of Nikol'skii classes $\mathbb{N}_p(\beta,L)$, see Definition~\ref{def:nikolskii-class},
where $\beta,L>0$ and $1\leq p<\infty$. It is important to emphasize that these classes allow estimating functions of inhomogeneous smoothness, i.e. those which can be very regular on some parts of the observation domain and rather irregular on the others.

\par The adaptation to the unknown parameters \( \ta^* \) and \( f \) can be formulated in terms of selection from a special family of  kernel estimators in the spirit of the Lepski and Goldenschluger-Lepski selection
  rules,  see \cite{Lep1990, KLP2001, GoldLep2009}. However, the proposed here procedure is quite different from
  the aforementioned ones, and it allows us to solve the problem of minimax adaptive estimation under the \( L_r \) losses over a collection of Nikol'skii classes.

 \par
In Section \ref{subsec:oracle-approach} we explain the proposed
selection rule and give the oracle inequalities.
Section \ref{subsec:adaptive-estimation} is devoted to the
application of these results to minimax adaptive estimation. 

\section{Oracle approach}
\label{subsec:oracle-approach}


Let $\K:\mathbb{R}\to\mathbb{R}$ be a function (kernel) satisfying $\int \K=1$. With any $\K$, any $z\in\mathbb{R}$, $h\in(0,1]$ and $f\in\mathbb{F}_M$ we associate \( \Delta_{\K,f}(h,z)=\sup_{\delta\leq h}\left|\delta^{-1}\int\K([u-z]/\delta)(f(u)-f(z))\dd u\right| \), a monotonous approximation error of the kernel smoother \( \delta^{-1} \int \K([u-z]/\delta)f(u)\dd u \). In particular, if the function $f$ is uniformly continuous, then\( \Delta_{\K,f}(h,z)\to 0 \) as \( h\to 0 \).
We will assume that the kernels are compactly supported symmetric Lipschitz functions.


\medskip

\par \noindent {\bf \emph{Oracle estimator.}}
For any $y\in\mathbb{R}$ define \( \overline{\Delta}_{\K,f}(h,y)=\sup_{a>0}(2a)^{-1}\int_{y-a}^{y+a}\Delta_{\K,f}(h,z)\dd z \), and $\Delta^{*}_{\K,f}(h,\cdot):=\max\left\{\overline{\Delta}_{\K,f}(h,\cdot),\Delta_{\K,f}(h,\cdot)\right\}$.
Define the oracle bandwidth:
\begin{equation}
\label{def_oracle bandwidth}
 \text{for any  } y \in \RR \quad
  h^*_{\K,f}(y) \eqdef
    \sup\big\{h\in [\eps^2, 1]:\; \sqrt{h}\;\Delta^{*}_{\K,f}(h,y)\leq
      \| \K \|_{\infty}\eps \sqrt{\ln (1/\eps)}\big\}.
\end{equation}
In what follows we assume that $ \eps\leq \exp{\big\{ -\max [1, (2M\|\K\|_1 \|\K\|_\infty^{-1})^2]\big\}}$. This assumption provides the well-defined $ h^*_{\K,f}$ and can be relaxed in several ways.

The quantity similar to the defined above \( h^*_{\K,f} \) first appeared in \cite{LepMamSpok97}
in the context of the estimating univariate functions possessing inhomogeneous smoothness.
Some years later this approach has been developed in \cite{KLP2001} and  \cite{GoldLep2009}
for multivariate function estimation.
In these papers, the interested reader can find a detailed discussion of the oracle approach.

\par For any $ (\theta,h)\in\SS\times[\eps^2, 1]$ define the matrix
$
E_{(\theta,h)}=\left(
\begin{array}{ll}
h^{-1}\theta_1\quad &h^{-1}\theta_2
\\
-\theta_2\quad &\;\theta_1
\end{array}
\right)
$, $\det(E_{(\ta, h)})=h^{-1}$,
and consider the family of kernel estimators
\begin{equation*}
\label{family of estimators}
\cc F
=
\Big\{ \hat F_{(\ta, h)} (\cdot)
    =
    \det(E_{(\ta, h)})\int K\big(E_{(\ta, h)}(t- \cdot)\big) Y_{\eps}(\dd t),\;\; (\theta,h)\in\SS\times[\eps^2, 1]\Big\}.
\end{equation*}
Here $K(u,v)=\K(u)\K(v)$, where $\K$ obeys the above conditions. Note that
\begin{equation}
\label{eq:distribution-of-kernel-estimator}
\hat F_{(\ta, h)} (\cdot)-\mathbb{E}^{\eps}_F\left[\hat{F}_{(\ta, h)} (\cdot)\right]\quad\sim\quad \mathcal{N}\left(0,\|\K\|^4_2\eps^{2} h^{-1}\right).
\end{equation}
The choice $\theta=\theta^*$ and $h=h^*:=h^*_{\K,f}(x^{T}\theta^*)$ leads to the {\it oracle }
(depending of \( F \)) estimator $\hat F_{(\theta^*, h^*)} (\cdot)$.
The meaning of this estimator is explained by the following result.
\begin{proposition}
\label{prop:risk-of-oracle-estimator}
For any \( (f,\theta^*)\in\mathbb{F}_M\times\SS \), \( r\ge 1 \) and \( \eps \) as above  we have
  \begin{equation*}
    \cc R_{r,x}^{(\eps)} \big(\hat{F}_{(\theta^*, h^*)}, F\big)\leq \mathfrak{c}_r
\| \K \|^{2}_{\infty}\eps \sqrt{\ln (1/\eps) / h^*_{\K,f}(x^{\T}\theta^*)}, \forall x\in[-1/2,1/2]^{2},
\end{equation*}
with $\mathfrak{c}_r=\left[\mathbb{E}\big(1+|\varsigma|\big)^{r}\right]^{1/r},\; \varsigma\sim \mathcal{N}(0,1)$.
\end{proposition}
This result means that the ``oracle'' knows the value of the index $\theta^*$ and the optimal, up to $\ln(1/\eps)$, trade-off $h^*$ between the approximation error determined by $\Delta^{*}_{\K,f}(h^{*},\cdot)$ and the stochastic error provided by the kernel estimator with the bandwidth $h^*$, cf.~(\ref{eq:distribution-of-kernel-estimator}).  That explains why the ``oracle'' chooses the ``estimator'' $\hat{F}_{(\theta^*, h^*)}$. Below we propose a ``real'' (based on the observation) estimator $\hat{F}(\cdot)$, which mimics the oracle. The construction of the estimator $\hat{F}(\cdot)$ is based on the data-driven selection from the family $\mathcal{F}$.

\smallskip

\par \noindent{\bf \emph{Selection rule.}} For any $\ta,\nu\in\SS$  and any $h\in[\eps^2,1]$ define the matrices
\begin{equation*}
\overline{E}_{(\theta,h)(\nu,h)}=\left(
\begin{array}{ll}
\frac{(\theta_1+\nu_1)}{2h(1+|\nu^\T\theta|)}\quad &\frac{(\theta_2+\nu_2)}{2h(1+|\nu^\T\theta|)}
\\*[2mm]
-\frac{(\theta_2+\nu_2)}{2(1+|\nu^\T\ta|)}\quad &\;\frac{(\theta_1+\nu_1)}{2(1+|\nu^\T\ta|)}
\end{array}
\right), \quad
E_{(\theta, h)(\nu,h)}=\left\{
\begin{array}{ll}
\overline{E}_{(\theta, h)(\nu,h)},\;\quad&\nu^\T\ta\geq 0;
\\*[2mm]
\overline{E}_{(-\theta, h)(\nu,h)},\quad&\nu^\T\ta< 0.
\end{array}
\right.
\end{equation*}
It is easy to check that
$
(4h)^{-1}\leq\det(E_{(\theta,h)(\nu,h)})\leq (2h)^{-1}.
$
The corresponding kernel estimator is defined by \( \hat F_{(\ta, h)(\nu, h)}(x)
    =    \det (E_{(\ta, h)(\nu,h)})\int K( E_{(\ta, h)(\nu,h)}(t-x))Y_{\eps}(\dd t) \).
For any $\eta\in (0,1]$ let \( \TH(\eta)= C(r, \K) \eps\sqrt{\eta^{-1}\ln (1/\eps)}\),
the constant \( C(r, \K) \) is given in \cite{LepSerd2011}, page~7.

Set
$\mathcal{H}_\eps=\big\{h_k=2^{-k},\; k = 0,1, \ldots\big\}\cap[\eps^2,1]$ and define for any \(\ta\in \SS \)
and  \(h \in \mathcal{H}_\eps \)
\begin{equation*}
R_{(\ta, h)}(x)
=
\sup_{\eta\in\mathcal{H}_\eps:\;  \eta\le h}
\Big\{
    \sup_{\nu\in\SS}
        \big|\hat F_{(\ta, \eta)(\nu, \eta)}(x) - \hat F_{(\nu, \eta)}(x) \big|   -\TH(\eta)
    \Big\}.
\end{equation*}
\noindent For any \( x \) introduce the random set
\( \cc P (x)=\{   (\ta, h)\in\SS\times\mathcal{H}_\eps : R_{(\ta, h)}(x) \le 0 \} \), and let
\( \tilde{h}=\max\{h:\, (\theta,h)\in\cc P (x)\} \) if  \( \cc P (x)\neq\emptyset \).
Note that there exists \( \vartheta\in\SS \) such that \( (\vartheta,\tilde{h})\in\mathcal{P}(x) \), since the set
\( \mathcal{H}_\eps \) is finite. Denote \( \hat{\Theta}:=\{\theta\in\SS:\,(\theta,\tilde{h})\in\cc P (x)\} \).
If \( \cc P (x)\neq\emptyset \), put \( \hat{\ta}= \ta \) such that \( \theta \in \hat{\Theta} \); otherwise
\(\hat{\ta}:= (1,0)^{\T} \). If \( \hat{\theta} \) is not unique, let us make any measurable choice. For instance, one can choose \( \hat{\theta}\in \hat{\Theta} \) with the smallest first coordinate.
Put as a final  estimator
\( \hat F (x)=\hat F_{(\hat \ta, \hat h)}(x)\), where
\begin{equation*}
    \hat{h}=\sup\left\{h\in\mathcal{H}_\eps:\;\;\left|\hat F_{(\hat{\ta},h)}(x)- \hat F_{(\hat{\ta}, \eta)}(x)\right|
                \leq\TH(\eta),\;\;\forall \eta\leq h,\;\eta\in\mathcal{H}_\eps\right\}.
\end{equation*}
\begin{theorem}\label{th:local-oracle-inequality}
{\bf \emph{Local and global oracle inequalities.}} For any $(f,\theta^*)\in\mathbb{F}_M\times\SS$, $r>0$
\begin{eqnarray*}
  \cc R_{r,x}^{(\eps)} \Big(\hat{F}_{(\hat{\theta}, \hat{h})}, F\Big)
  &\leq & C_{r,1}(\K)\sqrt{\frac{\|\K\|^4_\infty\eps^{2} \ln (1/\eps)}{h^*_{\K,f}(x^{T}\theta^*)}}
        +
            C_{r,2}(M,\K)\| \K \|^2_{\infty}\eps\sqrt{\ln (1/\eps)}, \, \forall x\in\left[-\frac{1}{2} ,\frac{1}{2} \right]^2; \\
  \cc R_r^{(\eps)} \Big(\hat{F}_{(\hat{\theta}, \hat{h})}, F\Big)
  &\leq &
  C_{r,1}(\K)\left\|\sqrt{\frac{\|\K\|^4_\infty\eps^{2}
 \ln (1/\eps)}{h^*_{\K,f}}}\right\|_r+C_{r,2}(M,\K)\| \K \|^2_{\infty}\eps\sqrt{\ln (1/\eps)}.
\end{eqnarray*}
The constants $C_{r,1}(\K)$ and $C_{r,2}(M,\K)$ are given in \cite{LepSerd2011}, page 11.
\end{theorem}
\section{Adaptive estimation}\label{subsec:adaptive-estimation}

In this section we apply the local oracle inequality given by the first assertion of Theorem \ref{th:local-oracle-inequality}
to the pointwise adaptive estimation over H\"older classes. Next, we use the global oracle inequality for adaptation over Nikol'skii classes.

\par For any \( a>0 \), denote by \( m_a  \), the maximal integer
strictly less than \( a \), and assume that there exists \( \bob>0 \) such that \( \int z^{j}K(z)\dd z=0,\;\;\forall j=1,\ldots,m_{\bob} \).

\bigskip

\par \noindent{\bf \emph{Pointwise adaptive estimation.}}
\begin{definition} \label{def:holder-class}
Let \( \beta>0 \) and \( L>0 \). A function \( g:\mathbb{R}\to\RR \) belongs to the H\"older class
\( \mathbb{H}(\beta,L) \), if \( g \) is \( m_\beta \)-times continuously
differentiable, \( \|g^{(m)}\|_\infty\le L,\;\forall m\leq m_\beta \),  and
$$
\left|g^{(m_\beta)}(t+h)-g^{(m_\beta)}(t)\right|\le L h^{\beta-m_\beta},\;\;\forall t \in\RR, \, h >0.
$$
\end{definition}

\par The aim is to estimate \( F(x) \) assuming that \( F\in\mathbb{F}(\bob):=\bigcup_{\beta\leq\bob}\bigcup_{L>0}\mathbb{F}_2(\beta,L) \), where
$$
\mathbb{F}_d(\beta,L)=\left\{F:\RR^d\to\RR\;|\; F(z)=f(z^{\T}\theta),\;f\in\mathbb{H}(\beta,L),\;\theta\in\SSd\right\},\, d\geq 2.
$$
Note that \( \bob \) can be an arbitrary number, but it must be chosen {\it a priory}.
\begin{theorem}
\label{th:pointwise-adaptation}
Let $\bob>0$ be fixed and the assumptions on the kernels hold. Then, for any
$\beta\leq\bob$, $L>0$, $x\in[-1/2,1/2]^{2}$, with $\psi_\eps(\beta,L)=L^{1/(2\beta+1)}\left(\eps\sqrt{\ln(1/\eps)}\right)^{2\beta/(2\beta+1)}$ we have
$$
\sup_{F\in\mathbb{F}_2(\beta,L)}\cc R_{r,x}^{(\eps)} \Big(\hat{F}_{(\hat{\theta}, \hat{h})}, F\Big)
\leq  \|\K \|^2_{\infty}\left[
C_{r,1}(\K)\psi_\eps(\beta,L)
+C_{r,2}(L,\K)\, \eps\sqrt{\ln (1/\eps)}\right].
$$
Moreover, for any $\beta,L>0$, $d\geq 2$  and  any $\eps>0$ small enough,
$$
\inf_{\tilde{F}}\sup_{F\in\mathbb{F}_d(\beta,L)}\cc R_{r,x}^{(\eps)} \Big(\tilde{F}, F\Big)
\geq \kappa \psi_\eps(\beta,L),
$$
where infimum is over all estimators. Here \( \kappa \) is a constant  independent of \( \eps \) and \( L \).
\end{theorem}
\noindent The estimator \( \hat{F}_{(\hat{\theta}, \hat{h})} \) is minimax adaptive with respect to
\( \{\mathbb{F}_d(\beta,L),\;\; \beta\le \bob,\; L>0\} \). It is surprising, since if the index is known, then
\( \mathbb{F}(\beta,L)=\mathbb{H}(\beta,L) \), and the problem can be reduced to the estimation of \( f \) at a point in the univariate Gaussian white noise model. As it is shown in \cite{Lep1990}
the optimally rate adaptive estimator over \( \{\mathbb{H}(\beta,L),\;\; \beta\leq \bob,\; L>0\} \) does not exist.

\bigskip

\par \noindent{\bf \emph{Adaptive estimation under the} \( \bb {L_r} \) \emph{losses.}}

\begin{definition}\label{def:nikolskii-class}
Let \( \beta>0 \), \( L>0 \), \( p\in [1,\infty) \). A function \( g:\mathbb{R}\to\RR \) belongs to the Nikol'skii class \( \mathbb{N}_p(\beta,L) \), if \( g \) is \( m_\beta \)-times continuously differentiable, 
\begin{equation*}
   \| g^{(m)} \|_p \le L,\; \forall m=1 \le m_\beta \quad \text{and} \quad
   \| g^{(m_\beta)}(\cdot+h)-g^{(m_\beta)}(\cdot) \|_p \le Lh^{\beta-m_\beta},\forall h>0.
\end{equation*}
We assume $\mathbb{N}_p(\beta,L)=\mathbb{H}(\beta,L)$ if $p=\infty$.
\end{definition}


\par Here the target of estimation is the function $F$  obeying the assumption \( F \in\mathbb{F}_p(\bob) \),
$ \mathbb{F}_p(\bob):=\bigcup_{\beta\leq\bob}\bigcup_{L>0}\mathbb{F}_{2,p}(\beta,L)$, where
\begin{eqnarray*}
&&\mathbb{F}_{d,p}(\beta,L)=\left\{F:\RR^d\to\RR\;|\; F(z)=f(z^{\T}\theta),\;f\in\mathbb{N}_p(\beta,L),\;\theta\in\SSd\right\}.
\end{eqnarray*}
\begin{theorem}
\label{th:global-adaptation}
Let $\bob>0$ be fixed and the above assumptions on the kernels hold.
Then
$$
\sup_{F\in\mathbb{F}_{2,p}(\beta,L)}\cc R_{r}^{(\eps)} \Big(\hat{F}_{(\hat{\theta}, \hat{h})}, F\Big)
\leq
\|\cc K\|^{2}_\infty \Big[ \kappa C_{r,1}(\K)\varphi_\eps(\beta,L,p)+C_{r,2}(L,\K)\eps\sqrt{\ln (1/\eps)}\Big],
$$
for any $L>0$, $p>1$, $p^{-1}<\beta\leq\bob$, and  $r\geq 1$. Here $\kappa$ is an absolute constant, and
$$
\varphi_\eps(\beta,L,p)=\left\{
\begin{array}{lll}
L^{1/(2\beta+1)}\left(\eps\sqrt{\ln(1/\eps)}\right)^{\frac{2\beta}{2\beta+1}},\quad& (2\beta+1)p>r;
\\
L^{1/(2\beta+1)}\left(\eps\sqrt{\ln(1/\eps)}\right)^{\frac{2\beta}{2\beta+1}}\big[\ln(1/\eps)\big]^{\frac{1}{r}},\quad& (2\beta+1)p=r;
\\
L^{\frac{1/2-1/r}{\beta-1/p+1/2}}
\left(\eps\sqrt{\ln (1/\eps)}\right)^{\frac{\beta-1/p+1/r}{\beta-1/p+1/2}},\quad& (2\beta+1)p< r.
\end{array}
\right.
$$
\end{theorem}
\noindent Note that $\mathbb{F}_{2,p}(\beta,L)\supset\mathbb{N}_{p}(\beta,L)$.
Indeed, the class  $\mathbb{N}_{p}(\beta,L)$ can be viewed as the class of functions $F$ satisfying $F(\cdot)=f(\ta^\T\cdot)$ with
$\ta=(1,0)^{\T}$. Then, the problem of estimating such (2-variate) functions  can be reduced  to the estimation
of univariate functions observed
in the one-dimensional GWN model. Thus, the rate of convergence for the latter problem, cf. \cite{DelyonJuditski1996, DJKP} and the references therein, is also the lower bound for the minimax risk  defined on
$\mathbb{F}_{2,p}(\beta,L)$.
Therefore the proposed estimator $\hat{F}_{(\hat{\theta}, \hat{h})}$
is optimally rate adaptive  whenever  $(2\beta+1)p<r$.
In the case $(2\beta+1)p\geq r$, we loose only a logarithmic factor with respect to the optimal rate,
and the construction of optimally rate adaptive estimator over a collection
$\big\{\mathbb{F}_{2,p}(\beta,L), \;\beta>0,\;L>0\big\}$
in this case remains an open problem.

\end{document}